\documentclass[11pt, oneside, secnum,leqno]{article}   	
\usepackage{geometry}                		
\usepackage{graphicx}				
\usepackage{amsmath,amssymb, amsbsy}
\newtheorem{thm}{Theorem}[section]
\newtheorem{crl}[thm]{Corollary}

\newtheorem{lmm}[thm]{Lemma}

\newtheorem{dfn}[thm]{Definition}
\newtheorem{exa}[thm]{Example}
\makeatletter
    
    \@addtoreset{equation}{section}
\makeatother

\title{Characterization of continuous endomorphisms in the space of entire functions of a given order}
\author{Takashi Aoki\footnote{Department of Mathematics, Kindai University, Higashi--Osaka 577-8502, Japan},
Ryuichi Ishimura\footnote{Department of Mathematics and Informatics Faculty of Science, Chiba University, Yayoicho, Chiba 263-8522,
Japan},
 Yasunori Okada\footnote{ Institute of Management and Information Technologies, Chiba University, Yayoicho, Chiba 263-8522,
Japan},
Daniele C. Struppa\footnote{Schmid College of Science and Technology, Chapman University, Orange 92866, CA, USA}
 \\ and Shofu Uchida\footnote{Graduate School of Science and Engineering, Kindai University, Higashi--Osaka 577-8502, Japan}
 }

\begin{document}
\maketitle

\section{Introduction}
The aim of this paper is to characterize continuous endomorphisms in the space of entire functions of
exponential type of order $p>0$. Let $A_p$ denote the space of entire functions of $n$ complex variables 
$z\in{\mathbb C}^n$ of order $p$ of normal type (cf. \cite{BG_book}, for example). 
We consider an endomorphism $F$ in the space, which is considered to be a DFS-space. 
We show that there is a unique linear differential operator $P$ of infinite order with coefficients in the space
which realizes $F$, that is, $Ff=Pf$ holds for any $f\in A_p$. The coefficients satisfy certain growth conditions
and conversely, if a formal differential operator of infinite order with coefficients in $A_p$ satisfy
these conditions, then it induces a continuous endomorphism. 

Characterization problem of linear continuous operators acting on a given function space has a long history.
J. Peetre (\cite{P}, \cite{P2}) proved that any morphism of sheaves of smooth functions becomes 
a differential operator locally of
finite order.  
Note that Peetre did not pose an assumption on the continuity.
In the analytic category, or in the ultradifferentiable category, operators of infinite order naturally appear
as continuous endomorphisms of function spaces. For example, any continuous sheaf
endomorphism (local operator) of the sheaf of holomorphic functions is expressed by a differential 
operators of infinite order satisfying certain estimates for coefficients. See \cite{SKK}, \cite{I1} for the analytic category, 
and \cite{Ko} for the ultradifferentiable category. For the spaces of entire functions of several variables,
characterization problem of continuous endomorphisms was firstly considered by the second author \cite{I2}. 
Some class of differential operators of infinite order was introduced and a part of characterization of
continuous endomorphisms was discussed. 
A similar class has been introduced in \cite{acss1}, \cite{acss2}  and
some applications to superoscillations (cf. references cited in these articles) have been given. 
We also note that several similar classes of infinite order differential operators with constant coefficients
or with polynomial coefficients were introduced and studied by \cite{vds} and our work is inspired by this thesis.
In this article, we give a complete answer to the characterization problem of the continuous 
endomorphisms in the space of entire functions of
a given order. 
A partial answer to this problem, for the case where $p>1$ and of normal type, is announced in \cite{aisu}.

The plan of this paper is as follows. 
In the second section, we recall some definitions of the spaces of entire functions and
give some refined estimates of the derivatives of entire functions. 
In the third section, we introduce a class of differential operators of infinite order and show that this is
exactly identified with the space of continuous endomorphisms of the space of entire functions
of a given order of normal type. In the fourth section, we give a similar result for the case of minimal type.
We give some examples and comments for generalization in the fifth section.

\section{The spaces of entire functions of a given order}
Let $p$ and $\tau$ be positive numbers. Let $A_{p,\tau}$ denote the set of all entire functions $f$ of
$n$ complex variables $z=(z_1,z_2,\dots,z_n)\in{\mathbb C}^n$ such that
\begin{equation}\label{ptaunorm}
\|f\|_{p,\tau}:=\sup_{z\in{\mathbb C}^n}|f(z)|\exp(-\tau |z|^p)<\infty.
\end{equation}
Here we set $|z|=\sqrt{|z_1|^2+|z_2|^2+\cdots|z_n|^2}$. 
We call $\|f\|_{p,\tau}$ the $(p,\tau)$-norm of $f$.
The set $A_{p,\tau}$ 
becomes a Banach space with the norm.
If $\tau<\tau'$, there is a natural inclusion mapping $A_{p,\tau}\hookrightarrow A_{p,\tau'}$. This mapping is a compact operator. Let $A_p$ (resp. $A_{p,0}$) denote the inductive limit (resp. the projective limit) 
of the family $\{A_{p,\tau}\}_{\tau>0}$ of the Banach spaces:
\begin{equation}\label{defAptau}
A_p:=\lim_{\xrightarrow[\tau>0]{}} A_{p,\tau} \quad (\mbox{resp.}\  A_{p,0}:=\lim_{\xleftarrow[\tau>0]{}} A_{p,\tau}).
\end{equation}
This becomes a DFS (resp. FS) space (cf. \cite{BG_book}). 
An entire function $f\in A_p$ (resp. $f\in A_{p,0}$) is said to be of order at most $p$ of normal type (resp.
of order at most $p$ of minimal type).
A linear operator $F:A_p\longrightarrow A_p$ 
is continuous if and only if
for any $\tau>0$ there exist $C>0$ and $\tau'>0$ for which $F(A_{p,\tau})\subset A_{p,\tau'}$ and
\[
\|Ff\|_{p,\tau'}\le C \|f\|_{p,\tau}
\] 
holds for any $f\in A_{p,\tau}$ (\cite[Chap 4, Part 1, 5, Corollary 1]{G}).

Similarly, a linear operator $F:A_{p,0}\longrightarrow A_{p.0}$ is continuous if and only if
for any  $\tau>0$, there exist $C>0$ and $\tau'>0$ for which
\[
\|Ff\|_{p,\tau}\le C \|f\|_{p,\tau'}
\]
holds for any $f\in A_{p,0}$.

To unify descriptions, we set $s_p=\max\{2^{p-1},1\}$. The following estimates for derivatives will be effectively used:
\begin{lmm}\label{firstlemma}
Let $f$ be an element of $A_{p,\tau}$. 
Then $\partial_z^\alpha f\in A_{p,s_p \tau}$ and we have
\begin{equation}\label{derivative0}
\|\partial_z^\alpha f\|_{p,s_p\tau}\le \frac{\alpha!}{|\alpha|^\frac{|\alpha|}{p}}(e\tau p)^{\frac{|\alpha|}{p}}
(2n^\frac{1}{2})^{|\alpha|}\|f\|_{p,\tau}.
\end{equation}
Here we set
\[
\partial_z^\alpha=\left(\frac{\partial}{\partial z_1}\right)^{\alpha_1}
\left(\frac{\partial}{\partial z_2}\right)^{\alpha_2}\cdots
\left(\frac{\partial}{\partial z_n}\right)^{\alpha_n}
\]
for a multi-index $\alpha=(\alpha_1,\alpha_2,\dots,\alpha_n)\in {\mathbb Z}_{\geq 0}^n$ and
$\alpha!:=\alpha_1!\alpha_2!\cdots\alpha_n!$, $|\alpha|:= \alpha_1+\alpha_2+\cdots+|\alpha_n|$.
\end{lmm}
{\it Proof} \ The Cauchy formula shows
\begin{equation}\label{derf}
|\partial_z^\alpha f(z)|\le \frac{\alpha!}{r^{|\alpha|}}\sup_{|\zeta_j|=r}|f(z+\zeta)|
\end{equation}
holds for any $r>0$, where we set $\zeta=(\zeta_1,\zeta_2,\dots,\zeta_n)$.
Since $f\in A_{p,\tau}$, we have
\begin{equation}\label{fA}
|f(z)|\le \|f\|_{p,\tau}\exp(\tau|z|^p)
\end{equation}
for all $z\in{\mathbb C}^n$.
Hence we have
\begin{equation}
|f(z+\zeta)|\le
\left\{
\begin{split}\label{fzzeta}
&\|f\|_{p,\tau}\exp(\tau(|z|^p+|\zeta|^p)) \quad \mbox{if}\ 0<p\le 1,\\
& \\
&\|f\|_{p,\tau}\exp(2^{p-1}\tau(|z|^p+|\zeta|^p)) \quad \mbox{if}\ 1<p.
\end{split}
\right.
\end{equation}
Here we have used the H\"older inequality for the case $p>1$.
Combining \eqref{derf} and \eqref{fzzeta}, we have
\begin{equation}
|\partial_z^\alpha f(z)|\le
\left\{
\begin{split}\label{derfbis}
&\frac{\alpha!}{r^{|\alpha|}}\|f\|_{p,\tau}\exp(\tau|z|^p)\exp(\tau n^\frac{p}{2} r^p) \quad \mbox{if}\  0<p\le 1,\\
& \\
&\frac{\alpha!}{r^{|\alpha|}}\|f\|_{p,\tau}\exp(2^{p-1}\tau|z|^p)\exp(2^{p-1}\tau n^\frac{p}{2} r^p) \quad \mbox{if}\ 1<p
\end{split}
\right.
\end{equation}
for all $r>0$. Taking the minimum of the right-hand sides of \eqref{derfbis}
with respect to $r$, we have \eqref{derivative0}.
\section{The case of normal type}
Let us consider a formal differential operator of infinite order with entire coefficients:
\begin{equation}\label{inf-ord}
P=\sum_\alpha a_\alpha(z)\partial_z^\alpha,
\end{equation}
where $\alpha$ runs over ${\mathbb Z}_{\ge 0}^n$. 
Let $p$ be a positive number and $q$ a real number satisfying $\displaystyle \frac{1}{p}+\frac{1}{q}=1$ if $p\neq 1$.
 If $p=1$, we set $\displaystyle \frac{1}{q}=0$.
For such an operator $P$, the following four conditions
are equivalent:
\begin{itemize}
\item[(I)] For any $\varepsilon>0$, there exist $B>0$ and $C>0$ for which
\begin{equation}\label{a-estim1}
|a_\alpha(z)|\le C\,\frac{\varepsilon^{|\alpha|}}{|\alpha|!^\frac{1}{q}}\exp(B|z|^p)
\end{equation}
holds for any $z\in{\mathbb C}^n$, $\alpha\in{\mathbb Z}_{\geq 0}^n$.  
\item[(II)] For any $\varepsilon>0$, there exist $B>0$ and $C>0$ for which
\begin{equation}\label{a-estim2}
|a_\alpha(z)|\le C\,\frac{\,\,|\alpha|^{\frac{|\alpha|}{p}}}{\alpha!}\,\varepsilon^{|\alpha|}\exp(B|z|^p)
\end{equation}
holds for any $z\in{\mathbb C}^n$, $\alpha\in{\mathbb Z}_{\geq 0}^n$. 
\item[(III)] For any $\varepsilon>0$, there exist $B>0$ and $C>0$ such that
\begin{equation}\label{a-estim3}
\|a_\alpha\|_{p,B}\le C\,\frac{\varepsilon^{|\alpha|}}{|\alpha|!^\frac{1}{q}}
\end{equation}
holds for any $\alpha\in{\mathbb Z}_{\geq 0}^n$. 
\item[(IV)] For any $\varepsilon>0$, there exist $B>0$ and $C>0$ such that
\begin{equation}\label{a-estim4}
\|a_\alpha\|_{p,B}\le \,C\,\frac{\,\,|\alpha|^{\frac{|\alpha|}{p}}}{\alpha!}\,\varepsilon^{|\alpha|}
\end{equation}
holds for any $\alpha\in{\mathbb Z}_{\geq 0}^n$.
\end{itemize}
Here and hereafter we use the standard convention for the multi-indices, that is, for $\alpha=(\alpha_1,\alpha_2,\dots,
\alpha_n)$ and $\beta=(\beta_1,\beta_2,\dots,\beta_n)\in{\mathbb Z}_{\ge 0}^n$, we set
$\alpha\pm\beta=(\alpha_1\pm\beta_1,\alpha_2\pm\beta_2,\dots,\alpha_n\pm\beta_n)$ and $z^\alpha=z_1^{\alpha_1}z_2^{\alpha_2}\cdots
z_n^{\alpha_n}$ for $z=(z_1,z_2,\dots,z_n)$, etc. We write $\beta\le \alpha$ if and only if $\beta_k\le \alpha_k$ for all $k=1,2,\dots,n$. 
The equivalence of the first two conditions follows from the following 
estimates which hold for any $\alpha\in{\mathbb Z}_{\ge 0}^n$:
\[
n^{-|\alpha|}e^{-\frac{|\alpha|}{p}}\frac{|\alpha|^\frac{|\alpha|}{p}}{\alpha!}
\le \frac{1}{|\alpha|!^\frac{1}{q}}\le \frac{|\alpha|^\frac{|\alpha|}{p}}{\alpha!}.
\]
The equivalence of the first and the third (resp. the second and the fourth) 
conditions comes from the definition of $(p,\tau)$-norm $\|\cdot\|_{p,\tau}$.
\begin{dfn}\label{Dp}. 
{\rm
The set of all formal differential operator $P$ of the form \eqref{inf-ord} satisfying one of the 
conditions (I)--(IV) is denoted by} $\boldsymbol{D}_p$.
\end{dfn}
\begin{thm}\label{mainthm1}
{\rm (i)} Suppose that $P\in\boldsymbol{D}_p$ has the form \eqref{inf-ord}. For an entire function $f\in A_p$, 
\[
Pf:=\sum_\alpha a_\alpha(z)\partial_z^\alpha f
\] 
converges and $Pf\in A_p$. Moreover, $f\mapsto Pf$ defines a linear continuous operator $P:A_p\longrightarrow A_p$.

\noindent {\rm (ii)} Let $F:A_p\longrightarrow A_p$ be a linear continuous operator. Then there is a unique 
$P\in\boldsymbol{D}_p$ such that $Ff=Pf$ holds for any $f\in A_p$. 
\end{thm}
{\bf Remark} We use the terminology ``linear continuous operator'' instead of ``continuous endomorphism''
which is included in the title of this paper. 
We use these two terms in the same meaning. The latter is shorter but the former seems to be
more familiar in functional analysis.

\

{\it Proof}\ \ Using Lemma \ref{firstlemma} 
and the estimation \eqref{a-estim4} of (IV) for $a_\alpha$, we have
for any $f\in A_{p,\tau}$,
\begin{equation}\label{Pfnorm}
\begin{split}
\|Pf\|_{p,B+s_p\tau} &\le \sum_{\alpha}\|a_\alpha\|_{p,B}\|\partial_z^\alpha f\|_{p,s_p\tau}\\
&\le C\sum_\alpha \frac{|\alpha|^\frac{|\alpha|}{p}}{\alpha!}\varepsilon^{|\alpha|}
\frac{\alpha!}{|\alpha|^\frac{|\alpha|}{p}}(e\tau p)^{\frac{|\alpha|}{p}}
(2n^\frac{1}{2})^{|\alpha|}\|f\|_{p,\tau}\\
&\le 2^{n-1}C\sum_{k=0}^\infty (e\tau p)^\frac{k}{p}(4\varepsilon n^\frac12)^k \|f\|_{p,\tau}.
\end{split}
\end{equation}
Here we have used the following two inequalities:
\[
\|fg\|_{p,\tau+\sigma}\le \|f\|_{p,\tau}\|g\|_{p,\sigma}\quad \mbox{for} \ \ f\in A_{p,\tau}, \ g\in A_{p,\sigma},
\]
\[
\sum_{|\alpha|=k} 1=\binom{n+k-1}{k}\le 2^{n+k-1}.
\]
If $\varepsilon$ is sufficiently small, the last sum in \eqref{Pfnorm} converges. For such an $\varepsilon>0$, 
we set $\tau'=B+s_p\tau$ and
\[
C'=\frac{2^{n-1}C}{1-4n^\frac12(e\tau p)^\frac{1}{p}\varepsilon}.
\]
Then we have
\[
\|Pf\|_{p,\tau'}\le C'\|f\|_{p,\tau}.
\]
Hence $Pf\in A_{p,\tau'}$ and we have obtained the continuity of $P:A_p\longrightarrow A_p$.

(ii) Let $F:A_p\longrightarrow A_p$ be a linear continuous operator. For any $\tau>0$, there exist
$C>0$ and $\tau'>0$ such that for any $f\in A_{p,\tau}$, we have $Ff\in A_{p,\tau'}$ and 
\begin{equation}\label{conti-est1}
\|F f\|_{p,\tau'}\le C\|f\|_{p,\tau}.
\end{equation}
Let us define a family of entire functions $\{a_\alpha(z)\}$ ($\alpha\in{\mathbb Z}_{\ge 0}^n$) by
\begin{equation}\label{a_alpha-def}
a_\alpha(z)=\sum_{\beta\le\alpha}\frac{(-1)^{|\alpha-\beta|} z^{\alpha-\beta}}{(\alpha-\beta)!\beta!}F z^\beta,
\end{equation}
where the convergence in $A_{p}$ shall be proved together with their estimates.
We define a formal differential operator $P$ of infinite order by
\begin{equation}\label{defP}
P=\sum_\alpha a_\alpha(z)\partial_z^\alpha.
\end{equation}
Firstly we show $P\in\boldsymbol{D}_p$. 
For any $\beta\in\mathbb Z_{\ge 0}^n$, we have
\begin{equation}\label{zbeta}
\|z^\beta\|_{p,\tau}\le\left(\frac{|\beta|}{e\tau p}\right)^\frac{|\beta|}{p}.
\end{equation}
In fact, $|z^\beta|\exp(-\tau|z|^p)$ is dominated by
\[
r^{|\beta|}\exp(-\tau r^{p})
\]
if $|z|=r$.
Taking the maximum of this function with respect to $r$, we obtain \eqref{zbeta}.
For any $\tau>0$, there exist $C>0$ and $\tau'>0$ such that
\begin{equation}\label{Fzbeta}
\|Fz^\beta\|_{p,\tau'}\le C\|z^\beta\|_{p,\tau}.
\end{equation}
Hence we have
\begin{equation}\label{aalphaestim}
\begin{split}
\|a_\alpha\|_{p,\tau+\tau'}&\le \sum_{\beta\le\alpha}\frac{\|z^{\alpha-\beta}\|_{p,\tau}}{(\alpha-\beta)!\beta!}
\|Fz^\beta\|_{p,\tau'}\\
&\le C\sum_{\beta\le\alpha}\frac{1}{(\alpha-\beta)!\beta!}
\left(\frac{|\alpha-\beta|}{e\tau p}\right)^\frac{|\alpha-\beta|}{p}\left(\frac{|\beta|}{e\tau p}\right)^\frac{|\beta|}{p}\\
&\le C\sum_{\beta\le\alpha}\frac{\alpha!}{(\alpha-\beta)!\beta!}\frac{|\alpha|^\frac{|\alpha|}{p}}{\alpha!}
\left(\frac{1}{e\tau p}\right)^\frac{|\alpha|}{p}\\
&\le C \frac{2^{|\alpha|}}{\alpha!}\left(\frac{|\alpha|}{e\tau p}\right)^\frac{|\alpha|}{p}.
\end{split}
\end{equation}
For a given $\varepsilon>0$, we can take $\tau>0$, $C>0$ and $\tau'>0$ so that \eqref{conti-est1} and
\[
\left(\frac{2^p}{e\tau p}\right)^\frac{1}{p}<\varepsilon
\]
hold.
If we set $B=\tau+\tau'$, we have
\[
\|a_\alpha\|_{p,B}\le C \frac{\ |\alpha|^\frac{|\alpha|}{p}}{\alpha!}\varepsilon^{|\alpha|}.
\]
for any $\alpha$. This implies $P\in\boldsymbol{D}_p$. 

To finish the proof, we compute $Pf$ for $f\in A_{p,\tau}$. We take the Taylor expansion of $f$:
\begin{equation}\label{taylor}
f(z)=\sum_{\mu} \frac{f_\mu}{\mu!}\,z^\mu \quad (f_{\mu}=\partial_{z}^{\mu} f(0)).
\end{equation}
For any $m\in\mathbb N$,  we can write
\begin{equation}\label{reminder}
f(z)-\sum_{|\mu|\le m-1}\frac{f_{\mu}}{\mu!}z^{\mu}=\frac{1}{(m-1)!}\int_{0}^{1}(1-t)^{m-1}\partial_{t}^{m}f(tz)dt.
\end{equation}
We set $f_{m,t}(z):=\partial_{t}^{m}f(tz)$.
We have
\begin{equation*}
\begin{split}
\|f_{m,t}\|_{p,s_p\tau+\tau'}&\le\sum_{|\alpha|=m}\frac{m!}{\alpha!}\|\partial_{z}^{\alpha}f\|_{p,s_p\tau}\|z^{\alpha}\|_{p,\tau'}\\
&\le m!2^{n-1}\left(\frac{4n^{\frac12}\tau^{\frac{1}{p}}}{\tau'^{\frac{1}{p}}}\right)^{m}\|f\|_{p,\tau}.
\end{split}
\end{equation*}
If $4n^{\frac12}\tau^{\frac{1}{p}}<\tau'^{\frac{1}{p}}$,
the $(p,s_p\tau+\tau')$-norm of the right-hand side of \eqref{reminder} converges to zero when $m\rightarrow\infty$.
This means that the Taylor expansion \eqref{taylor} converges in $A_{p,\tau''}$ for some $\tau''>0$. 
By the continuity of $F$, we have
\begin{equation}\label{Ff}
Ff=\sum_\mu \frac{f_\mu}{\mu!}\, Fz^\mu.
\end{equation}
On the other hand, by the definition of $a_\alpha$, we have
\begin{equation*}
\begin{split}
Pf&=\sum_\alpha a_\alpha(z)\sum_{\alpha\le\mu} \frac{f_\mu}{\mu!}\ \partial_z^\alpha z^\mu\\
&=\sum_{\beta\le\alpha\le\mu}\frac{(-z)^{\alpha-\beta}}{(\alpha-\beta)!\beta!}(Fz^\beta)\cdot 
\frac{f_\mu}{(\mu-\alpha)!}\, z^{\mu-\alpha}\\
&=\sum_\mu\frac{f_\mu}{\mu!}F z^\mu\\
&=Ff.
\end{split}
\end{equation*}
Here we have used an identity
\[
\sum_{\gamma\le\mu-\beta}\frac{(-1)^\gamma(\mu-\beta)!}{\gamma!(\mu-\beta-\gamma)!}=0
\]
which holds for any $\beta\le\mu$ satisfying $\beta\neq\mu$. This completes the proof.
\begin{crl}
The set $\boldsymbol{D}_p$  becomes a ring under natural addition and multiplication as differential operators.
\end{crl}
\section{The case of minimal type}
Let us consider a formal differential operator of infinite order with entire coefficients:
\begin{equation}\label{inf-ord0}
P=\sum_\alpha a_\alpha(z)\partial_z^\alpha,
\end{equation}
where $\alpha$ runs over ${\mathbb Z}_{\ge 0}^n$. 
Let $p$ be a positive number and $q$ a real number satisfying $\displaystyle \frac{1}{p}+\frac{1}{q}=1$ if $p\neq 1$.
 If $p=1$, we set $\displaystyle \frac{1}{q}=0$.
For such an operator $P$, the following four conditions
are equivalent:
\begin{itemize}
\item[$\rm(I)_0$] For any $\varepsilon>0$, there exist $B>0$ for which
\begin{equation}\label{a-estim01}
|a_\alpha(z)|\le \frac{B^{|\alpha|+1}}{|\alpha|!^\frac{1}{q}}\exp(\varepsilon|z|^p)
\end{equation}
holds for any $z\in{\mathbb C}^n$, $\alpha\in{\mathbb Z}_{\geq 0}^n$. 
\item[$\rm(II)_0$] For any $\varepsilon>0$, there exist $B>0$ for which
\begin{equation}\label{a-estim02}
|a_\alpha(z)|\le \frac{\,\,|\alpha|^{\frac{|\alpha|}{p}}}{\alpha!}\,B^{|\alpha|+1}\exp(\varepsilon|z|^p)
\end{equation}
holds for any $z\in{\mathbb C}^n$, $\alpha\in{\mathbb Z}_{\geq 0}^n$. 
\item[$\rm(III)_0$] For any $\varepsilon>0$, there exist $B>0$ such that
\begin{equation}\label{a-estim03}
\|a_\alpha\|_{p,\varepsilon}\le \frac{B^{|\alpha|+1}}{|\alpha|!^\frac{1}{q}}
\end{equation}
holds for any $\alpha\in{\mathbb Z}_{\geq 0}^n$. 
\item[$\rm(IV)_0$] For any $\varepsilon>0$, there exist $B>0$ such that
\begin{equation}\label{a-estim04}
\|a_\alpha\|_{p,\varepsilon}\le \,\frac{\,\,|\alpha|^{\frac{|\alpha|}{p}}}{\alpha!}\,B^{|\alpha|+1}
\end{equation}
holds for any $\alpha\in{\mathbb Z}_{\geq 0}^n$.
\end{itemize}
\begin{dfn}\label{Dp0}
{\rm
The set of all formal differential operator $P$ of the form \eqref{inf-ord0} satisfying one of the 
conditions $\rm(I)_0$--$\rm(IV)_0$ is denoted by} $\boldsymbol{D}_{p,0}$.
\end{dfn}
\begin{thm}\label{mainthm2}
{\rm (i)} Suppose that $P\in\boldsymbol{D}_{p,0}$ has the form \eqref{inf-ord0}. For an entire function $f\in A_{p,0}$, 
\[
Pf:=\sum_\alpha a_\alpha(z)\partial_z^\alpha f
\] 
converges and $Pf\in A_{p,0}$. Moreover, $f\mapsto Pf$ defines a linear continuous operator $P:A_{p,0}\longrightarrow A_{p,0}$.

\noindent {\rm (ii)} Let $F:A_{p,0}\longrightarrow A_{p,0}$ be a linear continuous operator. Then there is a unique 
$P\in\boldsymbol{D}_{p,0}$ such that $Ff=Pf$ holds for any $f\in A_{p,0}$. 
\end{thm}

\

{\it Proof}\ (i) We assume condition $\rm(IV)_0$ for $a_\alpha$. Similar 
computations in \eqref{Pfnorm} yield
\begin{equation}\label{Pfnorm0}
\begin{split}
\|Pf\|_{p,\varepsilon+s_p\tau}&\le\sum_{\alpha}\|a_\alpha\|_{p,\varepsilon}\|\partial_z^\alpha f\|_{p,s_p\tau}\\
&\le\sum_\alpha\frac{|\alpha|^\frac{|\alpha|}{p}}{\alpha!}B^{|\alpha|+1}
\frac{\alpha!}{|\alpha|^\frac{|\alpha|}{p}}(e\tau p)^\frac{|\alpha|}{p}(2n^\frac{1}{2})^{|\alpha|}\|f\|_{p,\tau}\\
&\le 2^{n-1}B\sum_{k=0}^\infty (e\tau p)^\frac{k}{p}(4Bn^\frac12)^k\|f\|_{p,\tau}
\end{split}
\end{equation}
for any $\varepsilon$, $\tau>0$, $f\in A_{p,0}$. For any $\sigma>0$, we can choose $\varepsilon>0$ and 
$\tau>0$ so that $\varepsilon+s_p\tau\le \sigma$ and $4(e\tau p)^\frac{1}{p}Bn^\frac12<1$. Thus there exist
$C>0$ and $\tau>0$ such that
\[
\|Pf\|_{p,\sigma}\le C\|f\|_{p,\tau}.
\]
This means $P:A_{p,0}\rightarrow A_{p,0}$ is continuous. 

(ii) For a given endomorphism $F:A_{p,0}\rightarrow A_{p,0}$, we construct 
$P=\sum_\alpha a_\alpha(z)\partial_z^\alpha$
in the same way as in Theorem \ref{mainthm1}, (ii), namely,
\[
a_\alpha(z)=\sum_{\beta\le\alpha}\frac{(-1)^{|\alpha-\beta|} z^{\alpha-\beta}}{(\alpha-\beta)!\beta!}F z^\beta.
\]
The right-hand side converges in $A_{p,0}$ for any $\alpha$. For any $\sigma>0$, there exist $C>0$ and 
$\sigma'>0$ so that $\|Fz^\beta\|_{p,\sigma}\le C\|z^\beta\|_{p,\sigma'}$ holds for all $\beta$.
Similar estimates to \eqref{aalphaestim} work and we have
\begin{equation}\label{aalphaestim0}
\begin{split}
\|a_\alpha\|_{p,\sigma+\tau}&\le \sum_{\beta\le\alpha}\frac{\|z^{\alpha-\beta}\|_{p,\tau}}{(\alpha-\beta)!\beta!}
\|Fz^\beta\|_{p,\sigma}\\
&\le C\sum_{\beta\le \alpha}\frac{1}{(\alpha-\beta)!\beta!}
\left(\frac{|\alpha-\beta|}{e\tau p}\right)^\frac{|\alpha-\beta|}{p}
\left(\frac{|\beta|}{e\sigma'p}\right)^\frac{|\beta|}{p}.
\end{split}
\end{equation}
For any $\varepsilon>0$, we choose $\sigma$, $\tau>0$ so that
$\sigma+\tau<\varepsilon$ holds. We can take $\tau$ as $\tau<\sigma'$. Hence 
$\|a_\alpha\|_{p,\varepsilon}\le\|a_\alpha\|_{p,\sigma+\tau}$ is
dominated by
\[
C \frac{2^{|\alpha|}}{\alpha!}\left(\frac{|\alpha|}{e\tau p}\right)^\frac{|\alpha|}{p}.
\]
If we set
\[
B=\max\left\{C, 2\left(\frac{1}{e\tau p}\right)^\frac{1}{p}\right\},
\]
we have
\[
\|a_\alpha\|_{p,\varepsilon}\le \frac{|\alpha|^\frac{|\alpha|}{p}}{\alpha!}B^{|\alpha|+1}.
\]
This implies $P\in{\boldsymbol D}_{p,0}$.
Convergence of the Taylor expansion of $f$ in $A_{p,0}$ can be proved similarly to Theorem \ref{mainthm1}, (ii).
Hence $Ff=Pf$ follows.
\begin{crl}
The set $\boldsymbol{D}_{p,0}$  becomes a ring under natural addition and multiplication as differential operators.
\end{crl}
\section{Examples and comments}
We give some examples of operators in $\boldsymbol{D}_{p}$ or $\boldsymbol{D}_{p,0}$.
Note that any differential operators of finite order with coefficients in $A_{p}$ (resp. $A_{p,0}$) belong to $\boldsymbol{D}_{p}$
(resp. $\boldsymbol{D}_{p,0}$) for any $p>0$. 
\begin{exa}
{\rm For $a=(a_{1},a_{2},\dots, a_{n})\in{\mathbb C}^{n}$, the translation operator $\displaystyle P=\sum_{\alpha} \frac{a^\alpha}{\alpha!}\partial_{z}^{\alpha}$ belongs to $\boldsymbol{D}_{p}$  and to $\boldsymbol{D}_{p,0}$ for any $p>0$. The action of $P$ is given as
$Pf(z)=f(z+a)$ for $f\in A_{p}$ or $f\in A_{p,0}$}
\end{exa}
We note that in \cite{acss1}, \cite{acss2}, a class of differential operators of 
infinite order was introduced, which was denoted by ${\mathcal D}_{p,0}$.
This is a subring of $\boldsymbol{D}_{p}$ given in Section 3 and the subscript ``$0$'' in this notation is not relevant to that of $\boldsymbol{D}_{p,0}$ given in Section 4.
\begin{exa}
{\rm Let $\sigma$ be a positive number and set
\[
P = \sum_{\alpha}\frac{\sigma^{|\alpha|}z^{\alpha}}{\alpha!}
    \partial_z^{\alpha}.
\]
This is a dilation operator, namely,
\[
Pf(z)=f((1+\sigma)z).
\]
It follows from \eqref{zbeta} that
this operator belongs to $\boldsymbol{D}_{p}$ for all $p>0$. 
We can also see that $P\notin\mathcal{D}_{p,0}$.}
\end{exa}
\begin{exa}{\rm
Let us consider the following initial value problem in ${\mathbb C}^2_{t,z}$:
\begin{equation}\label{sch}
\left\{
\begin{split}
&i\frac{\partial\psi}{\partial t}=\left(-\frac{1}{2}\frac{\partial^2}{\partial z^2}+z\right)\psi,\\
&\psi(0,z)=\phi(z).
\end{split}
\right.
\end{equation}
This is formally solved as
\begin{equation}
\psi(t,z)=\exp\left(\frac{t}{i}\left(-\frac{1}{2}\frac{\partial^2}{\partial z^2}+z\right)\right)\phi(z).
\end{equation}
The operator appeared in the right-hand side can be written in the form
\begin{equation}
P:=\exp\left(-i\left(tz+\frac{t^3}{6}\right)\right)\exp\left(\frac{it}{2}\frac{\partial^2}{\partial z^2}+\frac{t^2}{2}
\frac{\partial}{\partial z}\right).
\end{equation}
If $p<2$, this in an element of $\boldsymbol{D}_{p}$ on ${\mathbb C}_z$ for any $t$. Hence for every
$\phi\in A_p$, \eqref{sch} has a unique solution $\psi=P\phi$ which belongs to $A_p$ for any $t$.

}
\end{exa} 
{\bf Comments for generalization} 
In this article, we considered the spaces $A_p$ or $A_{p, 0}$ for a constant order $p$. The proofs of the results obtained in Sections 3 and 4 are in fact valid for the case of the spaces of entire functions of exponential type with respect to a given proximate order $\rho(r)$
introduced by Valiron \cite{V} in one variable case (for such spaces in several variable case, see for example \cite{L-G}) and we can give
the complete answer to the corresponding characterization problem of the continuous endomorphisms.
We remark that in \cite{I2}, one of the authors considered this problem but gave only a partial answer which is incomplete comparing with the results of present article.

\

\noindent{\bf Acknowledgments}

This work was supported by JSPS KAKENHI Grant Numbers 26400126,  16K05170 and 18K03385.

\end{document}